\documentclass{amsart}
\usepackage[left=25mm, right=25mm, top=1in, bottom=1in]{geometry}
\usepackage{amsmath,bm,amssymb}
\usepackage{scrextend} 
\allowdisplaybreaks
\usepackage[normalem]{ulem}
\usepackage[dvipsnames]{xcolor}
\usepackage{setspace}
\usepackage{amsbsy}
\usepackage{hyperref}
\usepackage{graphicx}
\usepackage[overload]{empheq}

\theoremstyle{plain}
\newtheorem{theorem}{Theorem}

\newtheorem{definition}[theorem]{Definition}

\newcommand{\del}{\partial}
\newcommand{\diver}{\operatorname{div}}

\hypersetup{
    colorlinks=true,
    linkcolor=blue,
    filecolor=blue,      
    urlcolor=blue,
}

\numberwithin{equation}{section}

\subjclass[2020]{35D30, 35Q35, 76N10, 76R50, 76T30.}
        \keywords{Multicomponent, Maxwell--Stefan, Cross-diffusion}  

\begin{document}

\title{Global existence of weak solutions for the \\
Maxwell--Stefan system in the whole space}

\author[S. Georgiadis]{Stefanos Georgiadis}
        \address[S. Georgiadis]{Division of Science and Mubadala Arabian Center for Climate and Environmental Science\\ New York University
        Abu Dhabi \\ United Arab Emirates}
        \email[]{\href{seg577@nyu.edu}{seg577@nyu.edu}}

\begin{abstract}
    We prove the global-in-time existence of weak solutions to the isothermal Maxwell--Stefan system on the whole space $\mathbb R^3$. The main difficulty is that, unlike in bounded domains, the concentrations generally have infinite mass and the standard mixing entropy is not finite. We therefore work with the relative entropy with respect to a strictly positive constant equilibrium state. The proof proceeds by solving the problem on balls $B_R$ with no-flux boundary conditions, using the bounded-domain entropy theory, and deriving estimates independent of $R$. These estimates yield uniform control of the relative entropy, the gradients $\nabla\sqrt{c_i}$, and the fluxes $J_i$. Passing to the limit $R\to\infty$ is achieved by local compactness and a diagonal argument. The resulting weak solution satisfies the Maxwell--Stefan system in the sense of distributions and obeys a global relative entropy inequality.
\end{abstract}

\maketitle



\section{Introduction}

The Maxwell--Stefan system is a fundamental model for multicomponent mass transport in fluid mixtures. It is derived from a balance between thermodynamic driving forces and interspecies friction, rather than from a direct proportionality between fluxes and concentration gradients as in Fickian diffusion. The dependence of each species flux on the state and motion of the other species leads to a strongly coupled cross-diffusion system, in which the evolution of one component is influenced by the concentration gradients of all components. This intrinsic coupling between species makes the Maxwell--Stefan approach particularly suitable for modeling transport processes in gases, liquids, porous media, membrane separations, and electrochemical devices.

Consider the isothermal Maxwell--Stefan system 
\begin{equation}\label{mass}
    \partial_t c_i+\diver J_i=0, \qquad i\in\{1,\dots,n+1\},   
\end{equation}
where $c=(c_1,\dots,c_{n+1})$ is the vector of concentrations and $J=(J_1,\dots,J_{n+1})$ the vector of mass fluxes, which are determined by solving uniquely the constrained linear system
\begin{align}
    - \sum_{j\not=i}^{n+1} \frac{c_j J_i-J_j c_i}{D_{ij}} &= c_i \nabla \mu_i, \quad i\in\{1,\dots,n+1\} \label{linearsystem} \\
    \sum_{i=1}^{n+1} J_i &= 0 \label{constraint},
\end{align}
where $D_{ij}$ are positive symmetric interaction coefficients and $\mu_i$ denotes the chemical potential of the $i$-th species, given by the variational derivative of the associated entropy (or free energy) functional $H$, namely
\[
\mu_i = \frac{\delta H}{\delta c_i}.
\]
The concentrations $c_i$ satisfy the bounds $0 \leq c_i \leq 1$, as well as the normalization constraint
\begin{equation} \label{normal}
    \sum_{i=1}^{n+1} c_i=1.
\end{equation}

The mathematical analysis of the Maxwell--Stefan system has attracted considerable attention in recent years. The existence of local-in-time strong solutions has been established in \cite{GM98} for the Cauchy problem and in \cite{B11} for bounded domains, while weak solutions were shown to exist globally in time in \cite{JS13} on a bounded subset of $\mathbb{R}^3$. For results on the nonisothermal case we refer to \cite{HJ21, GJ24} for weak solutions and \cite{HMPW17, HS18, MT15} for strong solutions.

The critical estimates for the existence theory of weak solutions come from the idea of associating an entropy functional to the system, which plays the role of a Lyapunov functional. In the case of a bounded domain $\Omega$, the entropy is the standard Boltzmann/mixing entropy, namely
\begin{equation}\label{mixingentr}
    H_0[c(t)]=\int_\Omega h_0(c)\,dx=\int_\Omega \sum_{i=1}^{n+1} c_i(\log{c_i}-1) \,dx.    
\end{equation}
The objective of the present work is to show the existence of global-in-time weak solutions on the whole space $\mathbb{R}^3$. Contrary to the bounded domain case, one cannot use the same entropy because $c_i$ do not have finite mass \cite{D8}. Instead, we fix a constant equilibrium state $\bar{c}=(\bar c_1,\dots,\bar c_{n+1})$, $\bar{c}_i>0$, $\sum_{i=1}^{n+1}\bar{c}_i=1$, and use the relative entropy
\begin{equation} \label{relentr}
    H[c|\bar c]=\int_{\mathbb{R}^3}h(c|\bar c)\,dx=\int_{\mathbb{R}^3}\sum_{i=1}^{n+1}c_i\log\frac{c_i}{\bar c_i}\,dx.
\end{equation}

The use of the relative entropy does not change the underlying Maxwell--Stefan dynamics. Indeed, the relative entropy density differs from the standard mixing entropy density only by a constant and an affine term in the concentrations. More precisely, since
\[
h_0(c)=\sum_{i=1}^{n+1} c_i(\log c_i-1),
\]
we have
\[
h(c|\bar c) = \sum_{i=1}^{n+1}c_i\log\frac{c_i}{\bar c_i} = h_0(c)+1-\sum_{i=1}^{n+1}c_i\log \bar c_i,
\]
where we used the constraint $\sum_i c_i=1$. Hence the corresponding chemical potentials differ only by constants, and therefore their gradients coincide:
\[
\nabla\log\frac{c_i}{\bar c_i}=\nabla\log c_i.
\]
Thus the relative entropy is a natural Lyapunov functional for the whole-space problem, while leading to the same Maxwell--Stefan flux relations.

The main difficulty in passing from bounded domains to the whole space is twofold. First, the concentrations are not expected to have finite mass on $\mathbb R^3$. In fact, if $c_i$ approaches a strictly positive constant state $\bar c_i$ at spatial infinity, then
\[
\int_{\mathbb R^3} c_i\,dx=\infty.
\]
Consequently, the standard entropy \eqref{mixingentr} is not the appropriate quantity on $\mathbb R^3$. Instead, a finite-entropy condition is imposed on the perturbation from equilibrium through the relative entropy \eqref{relentr}. Second, the compact embedding $H^1(\Omega)\hookrightarrow L^2(\Omega)$, which is available on bounded domains, has no direct analogue on $\mathbb R^3$. Hence one cannot pass to the limit by using global compactness. The proof must instead rely on local compactness and a careful control of the behavior at infinity.

The purpose of this paper is to prove the global-in-time existence of weak solutions to the isothermal Maxwell--Stefan system on $\mathbb R^3$ for initial data with finite relative entropy with respect to a strictly positive constant equilibrium state. More precisely, assuming
\[
0\leq c_i^0\leq 1,\qquad \sum_{i=1}^{n+1}c_i^0=1,
\]
and
\[
H[c^0|\bar c]<\infty,
\]
we prove that, for every $T>0$, there exists a weak solution $(c,J)$ on $(0,T)\times\mathbb R^3$ satisfying
\[
0\leq c_i\leq 1,\qquad \sum_{i=1}^{n+1}c_i=1,\qquad \sum_{i=1}^{n+1}J_i=0,
\]
as well as the estimates
\[
c_i-\bar c_i\in L^\infty(0,T;L^2(\mathbb R^3)),
\]
\[
\nabla\sqrt{c_i}\in L^2((0,T)\times\mathbb R^3;\mathbb{R}^3),
\]
and
\[
J_i\in L^2((0,T)\times\mathbb R^3;\mathbb R^3).
\]
Moreover, the solution satisfies the relative entropy inequality
\[
H[c(t)|\bar c]+C\int_0^t\int_{\mathbb R^3} \sum_{i=1}^{n+1}|\nabla\sqrt{c_i}|^2\,dx\,ds \leq H[c^0|\bar c]
\]
for almost every $t>0$, where $C>0$ depends only on the Maxwell--Stefan diffusivities.

The proof follows the entropy method of Jüngel and Stelzer \cite{JS13} on bounded domains, but with several modifications required by the whole-space setting. We first solve the problem on balls $B_R\subset\mathbb R^3$ with homogeneous no-flux boundary conditions. The bounded-domain theory gives weak solutions on each $B_R$. Since the no-flux condition conserves the mass of each species on $B_R$, the standard entropy inequality can be rewritten as a relative entropy inequality with respect to the fixed constant state $\bar c$. This yields
\[
H_{B_R}[c^R(t)|\bar c] + \int_0^t \mathcal D_{B_R}[c^R,J^R]\,ds \leq
H_{B_R}[c^{0,R}|\bar c] \leq H[c^0|\bar c],
\]
This uniform estimate is the key point of the argument.

The entropy dissipation controls both the gradients of the square roots of the concentrations and the fluxes. More precisely, using the Maxwell--Stefan relations and the constraint $\sum_i J_i=0$, one obtains bounds of the form
\[
\int_0^T\int_{B_R} \sum_{i=1}^{n+1}|\nabla\sqrt{c_i^R}|^2\,dx\,dt \leq C H[c^0|\bar c], 
\]
and
\[
\int_0^T\int_{B_R} \sum_{i=1}^{n+1}|J_i^R|^2\,dx\,dt \leq C  H[c^0|\bar c],
\]
where $C$ is independent of $R$. These estimates are purely entropic and algebraic; in particular, no Poincaré inequality on $B_R$ is used. This is essential, since the Poincaré constant would diverge as $R\to\infty$.

The passage to the limit $R\to\infty$ is obtained by local compactness. On every fixed ball $B_K$, the estimates imply
\[
c_i^R-\bar c_i \quad\text{bounded in}\quad L^\infty(0,T;L^2(B_K)),
\]
\[
c_i^R \quad\text{bounded in}\quad L^2(0,T;H^1(B_K)),
\]
and
\[
\partial_t c_i^R = -\operatorname{div}J_i^R \quad\text{bounded in}\quad L^2(0,T;W^{-1,2}(B_K)).
\]
The Aubin--Lions lemma then gives strong convergence of $c_i^R$ in $L^2((0,T)\times B_K)$. A diagonal argument over $K$ yields strong  local convergence on $\mathbb R^3\times(0,T)$, which is sufficient to pass to the nonlinear Maxwell--Stefan relations. Finally, lower semicontinuity of the relative entropy and of the dissipation gives the entropy inequality on the whole space.

We emphasize that the far-field condition is not imposed pointwise. The condition
\[
H[c^0|\bar c]<\infty
\]
is the natural finite-entropy replacement for finite mass in the present setting. It implies, in particular, that 
\[
c_i^0-\bar c_i\in L^2(\mathbb R^3),
\]
and the entropy inequality propagates this information in time. Thus the solution remains close to the equilibrium state in the relative-entropy sense, although the individual concentrations themselves generally have infinite mass on $\mathbb R^3$.

Unlike in bounded domains, we do not rely on Poincaré-type inequalities and we do not obtain exponential convergence to equilibrium. The result is therefore an existence theorem for global finite-entropy weak solutions on the whole space, rather than a long-time relaxation result.


\section{Assumptions and main results}

We impose the following assumption:

\begin{itemize}
    \item[(A1)] We work on $\mathbb{R}^3 \times (0,T)$ for some $T>0$.
    \item[(A2)] There exists a constant equilibrium state $\bar c=(\bar c_i, \dots, \bar c_{n+1})$, with $\bar c_i>0$ and $\sum_{i=1}^{n+1}\bar c_i = 1$.
    \item[(A3)] The initial data $c(x,0)=(c_1^0(x),\dots,c_{n+1}^0(x))$ are chosen such that $0 \leq c_i^0 \leq 1$ and $\sum_{i=1}^{n+1} c_i^0=1$, with $H[c^0|\bar c] < +\infty$.
\end{itemize}

\begin{definition}[Weak solutions] \label{def}
    A weak solution to \eqref{mass}--\eqref{constraint} consists of functions $c_i:\mathbb{R}^3\times(0,T) \to [0,1]$, such that
    \[
    \sum_{i=1}^{n+1}c_i=1, \quad \sum_{i=1}^{n+1}J_i=0
    \]
    with the regularity
    \[
    (c_i - \bar{c}_i) \in L^\infty(0,T;L^2(\mathbb{R}^3)), \quad \nabla \sqrt{c_i} \in L^2(\mathbb{R}^3\times(0,T);\mathbb{R}^3) , \quad J_i \in L^2(\mathbb{R}^3\times(0,T);\mathbb{R}^3)
    \]
    for all $i \in\{1,\dots,n+1\}$, where $J_i:\mathbb{R}^3\times(0,T)\to\mathbb{R}^3$ satisfies the linear system \eqref{linearsystem} in the sense of distributions and the constraint \eqref{constraint} almost everywhere, such that
    \begin{align}
        \int_0^T \int_{\mathbb{R}^3} c_i \del_t\phi \,dxdt + \int_{\mathbb{R}^3}c_i^0(x) \phi(x,0) \,dx + \int_0^T \int_{\mathbb{R}^3} J_i \cdot \nabla\phi \, dxdt = 0, \label{weakmass}
    \end{align}
    for test functions $\phi \in C_c^\infty(\mathbb{R}^3\times[0,T))$.
\end{definition}

\begin{theorem}[Existence of weak solutions] \label{thmex}
    Let assumptions (A1)--(A3) hold. For every $T>0$, there exists a weak solution to \eqref{mass}--\eqref{constraint} according to Definition \ref{def}, which satisfies the entropy inequality 
    \[
        H[c(t)|\bar{c}] + C\int_0^t\int_{\mathbb{R}^3}\sum_{i=1}^{n+1}|\nabla\sqrt{c_i}|^2 \,dxdt \leq H[c^0|\bar{c}]
    \]
    for almost every $t\in(0,T)$.
\end{theorem}


\section{Proof of Theorem \ref{thmex}}

\subsection*{Preliminaries} The proof is based on exhaustion of $\mathbb{R}^3$ by balls. Let 
\[
B_R = \{x\in\mathbb{R}^3 : |x|< R \}.
\]
For fixed $R>0$, we consider the Maxwell--Stefan system on $B_R$, with no-flux boundary conditions
\[
J_i^R \cdot \nu=0 \quad \textnormal{on } \partial B_R
\]
and the initial conditions
\[
c^{0,R}=c^0(x)\Big|_{B_R}
\]
which satisfy
\[
0\leq c_i^{0,R}\leq 1, \quad \sum_i c_i^{0,R} = 1, \quad \textnormal{ a.e. in } B_R.
\]
Then
\begin{equation}\label{unifentr}
    H_{B_R}[c^{0,R}|\bar c]=\int_{B_R}h(c^{0,R}|\bar c) \,dx \leq \int_{\mathbb{R}^3}h(c^{0}|\bar c) \,dx=H[c^{0}|\bar c]
\end{equation}
i.e. all approximate initial entropies are bounded by the same quantity $H[c^{0}|\bar c]$, which is independent of $R$.

By the result in \cite{JS13} for bounded domains, we obtain, for each $R>0$, a global bounded weak solution $(c^R,J^R)$ on $B_R\times(0,T)$. Thus
\[
0\leq c_i^R\leq1, \quad \sum_i c_i^R=1, \quad \sum_iJ_i^R=0.
\]
Strictly speaking, the existence theorem of \cite{JS13} is based on the entropy functional \eqref{mixingentr}, but we are using the relative entropy \eqref{relentr}. The reason the relative entropy causes no conflict is that it differs from the usual entropy only by an affine term in $c$. Indeed,
\[
h(c|\bar c) = h_0(c) + 1 - \sum_ic_i\log{\bar c_i}.
\]
Because of the no-flux boundary condition on $B_R$, the mass of each species is conserved on $B_R$:
\[
\frac{d}{dt}\int_{B_R}c_i^R\,dx = -\int_{\partial B_R} J_i^R\cdot\nu \,dS = 0
\]
and so 
\[
\int_{B_R}c_i^R(t)\,dx=\int_{B_R}c_i^{0,R}\,dx.
\]
Therefore, due to $\bar c$ being constant, the affine correction to the entropy 
\[
|B_R| - \sum_i\log{\bar c_i}\int_{B_R}c_i^R(t)\,dx
\]
is constant in time. Hence, the entropy variables are merely shifted by constants and the Hessian of the entropy remains the same. Thus, the parabolic structure, the coercivity, and the entropy-variable method are unchanged. Consequently, the entropy inequality of \cite{JS13} is equivalent to the relative entropy inequality 
\[
H_{B_R}[c^R|\bar c] + \int_0^t D_R[c^R,J^R] \,dt \leq H_{B_R}[c^{0,R}|\bar c],
\]
where $D_R[c^R,J^R]$ is the dissipation, which reads
\[
D_R[c^R,J^R] = \frac{1}{2}\int_{B_R} \sum_{i=1}^{n+1}\sum_{j\not=i}\frac{c_i^Rc_j^R}{D_{ij}}|u_i^R-u_j^R|^2 \,dx
\]
with
\[
J_i^R=c_i^Ru_i^R, \quad \sum_ic_i^Ru_i^R=0.
\]
The velocities $u_i$ are not needed at the final level, but they are useful in order to derive the estimates.

Due to the uniform bound of the initial approximate entropies \eqref{unifentr}, we obtain the uniform-in-$R$ estimate
\begin{equation} \label{entrineq}
        H_{B_R}[c^R|\bar c] + \int_0^t D_R[c^R,J^R] \,dt \leq H[c^0|\bar c],
\end{equation}

\subsection*{Estimates for fixed $R$} Pinsker's inequality (sometimes also called Kullback-Leibler inequality), or simply the convexity of the relative entropy, implies
\[
|c-\bar c|^2 \leq 2h(c|\bar c).
\]
Hence
\[
\int_{B_R}|c^R-\bar c|^2 \,dx \leq \int_{B_R} 2h(c^R|\bar c) \,dx
\]
and the right-hand side is uniformly bounded by the initial relative entropy due to the relative entropy inequality \eqref{entrineq}. Taking the essential supremum over time, we find
\[
\|c^R-\bar c\|_{L^\infty(0,T;L^2(B_R))} \leq \sqrt{2 H[c^0|\bar c]}.
\]

Moreover, the dissipation provides gradient estimates. Indeed, using \eqref{linearsystem} and Cauchy-Schwarz, we find
\begin{align*}
    |\nabla\sqrt{c_i^R}|^2 & = \frac{|\nabla c_i^R|^2}{4 c_i} = \frac{c_i^R}{4} \left|\sum_{j\not=i} \frac{c_j^R}{D_{ij}}(u_i^R-u_j^R)\right|^2 \\
    & \leq \frac{c_i^R}{4} \left(\sum_{j\not=i} \frac{c_j^R}{D_{ij}}\right) \left(\sum_{j\not=i} \frac{c_j^R}{D_{ij}}|u_i^R-u_j^R|^2\right) \\
    & \leq C\sum_{j\not=i} \frac{c_i^Rc_j^R}{2D_{ij}}|u_i^R-u_j^R|^2
\end{align*}
and thus
\begin{equation} \label{grad}
    \int_0^T\int_{B_R}\sum_i|\nabla\sqrt{c_i^R}|^2\,dxdt \leq C\int_0^T D_R[c^R,J^R] \,dt \leq C H[c^0|\bar c]
\end{equation}
where the constant does not depend on $R$.

Since
\[
\nabla c_i^R = 2\sqrt{c_i^R} \nabla \sqrt{c_i^R}
\]
we also find
\[
\int_0^T\int_{B_R}\sum_i|\nabla c_i^R|^2\,dxdt \leq C H[c^0|\bar c]
\]
for a different constant independent of $R$.

Finally, the dissipation also controls the fluxes. Due to the constraint \eqref{constraint}, one has
\[
u_i^R = \sum_{j\not=i} c_j^R (u_i^R-u_j^R)
\]
and by Jensen's inequality
\[
|u_i^R|^2 \leq \sum_{j\not=i} c_j^R |u_i^R-u_j^R|^2.
\]
This suggests that
\[
|J_i^R|^2 = |c_i^Ru_i^R|^2 \leq \sum_{j\not=i} c_i^Rc_j^R |u_i^R-u_j^R|^2 \leq \max_{i\not=j} D_{ij} \sum_{j\not=i} \frac{c_i^Rc_j^R}{D_{ij}} |u_i^R-u_j^R|^2  
\]
and so
\[
\int_0^T\int_{B_R}\sum_i|J_i^R|^2\,dxdt \leq C H[c^0|\bar c].
\]
Again the constant is independent of $R$.

\subsection*{Extension to $\mathbb{R}^3$} We, now, extend $c^R$ by the equilibrium state $\bar c$ outside the ball $B_R$
\[
\tilde c^R=\begin{cases}
        c^R(x,t), & x\in B_R \\
        \bar c, & x\not\in B_R
\end{cases}
\]
and $J^R$ by zero
\[
\tilde J^R=\begin{cases}
        J^R(x,t), & x\in B_R \\
        0, & x\not\in B_R  
\end{cases}.
\]
From the previous estimates
\[
\tilde c_i^R - \bar c_i \textnormal{ is bounded in } L^\infty(0,T;L^2(\mathbb{R}^3))
\]
and
\[
\tilde J_i^R \textnormal{ is bounded in } L^2(\mathbb{R}^3\times(0,T);\mathbb{R}^3).
\]
Therefore, after extracting a subsequence,
\[
(\tilde c_i^R - \bar c_i) \overset{*}{\rightharpoonup} (c_i-\bar c_i) \textnormal{ weakly-$*$ in } L^\infty(0,T;L^2(\mathbb{R}^3))
\]
and
\[
\tilde J_i^R \rightharpoonup J_i \textnormal{ weakly in } L^2(\mathbb{R}^3\times(0,T);\mathbb{R}^3)
\]
i.e.
\[ 
(c_i-\bar c_i) \in L^\infty(0,T;L^2(\mathbb{R}^3)), \quad J_i \in L^2(\mathbb{R}^3\times(0,T);\mathbb{R}^3).
\]

\subsection*{Local compactness} Fix $K>0$. For $R>K$, the functions $c_i^R$ are defined on $B_K$. From the estimates above, we find
\[
c_i^R - \bar c_i \textnormal{ is bounded in } L^\infty(0,T;L^2(B_K))
\]
\[
\nabla c_i^R \textnormal{ is bounded in } L^2(B_K\times(0,T);\mathbb{R}^3).
\]
and
\[
J_i^R \textnormal{ is bounded in } L^2(B_K\times(0,T);\mathbb{R}^3).
\]
thus
\[
c_i^R \textnormal{ is bounded in } L^2(0,T;H^1(B_K))
\]
and
\[
\del_t c_i^R = -\diver J_i^R \textnormal{ is bounded in } L^2(0,T;W^{-1,2}(B_K)).
\]
By Aubin--Lions lemma
\[
c_i^R \to c_i \textnormal{ strongly in } L^2(B_K\times(0,T))
\]
and after a diagonal extraction $K=1,2,\dots$
\[
c_i^R \to c_i \textnormal{ strongly in } L^2_{\rm loc}(\mathbb{R}^3\times(0,T)).
\]
Passing to a further subsequence if necessary, we also conclude
\[
c_i^R \to c_i \textnormal{ a.e. locally in } \mathbb{R}^3\times(0,T).
\]
This allows us to pass to the limit in
\[
0\leq c_i^R\leq1, \quad \sum_i c_i^R=1
\]
and obtain 
\[
0\leq c_i \leq1, \quad \sum_i c_i=1.
\]
Also, since
\[
\sum_i J_i^R=0
\]
and $J_i^R\rightharpoonup J_i$ weakly in $L^2_{\rm loc}$, we conclude
\[
\sum_i J_i=0 \textnormal{ a.e. in } \mathbb{R}^3\times(0,T).
\]

\subsection*{Passage to the limit in the weak formulation}

Let $\phi \in C_c^\infty(\mathbb{R}^3\times[0,T))$ be a test function. Choose $K>0$, such that $\textnormal{supp}\phi \subset B_K\times[0,T)$. For any $R>K$, the weak formulation for a bounded domain reads
\[
\int_0^T \int_{B_R} c_i^R \del_t\phi \,dxdt + \int_{B_R}c_i^{0,R}(x) \phi(x,0) \,dx + \int_0^T \int_{B_R} J_i^R \cdot \nabla\phi \, dxdt = 0, 
\]
but because $\phi$ is supported in $B_K$, we may write
\[
\int_0^T \int_{\mathbb{R}^3} c_i^R \del_t\phi \,dxdt + \int_{\mathbb{R}^3}c_i^0(x) \phi(x,0) \,dx + \int_0^T \int_{\mathbb{R}^3} J_i^R \cdot \nabla\phi \, dxdt = 0.
\]
Using the convergence
\[
c_i^R \to c_i \textnormal{ strongly in } L^2(\textnormal{supp}\phi)
\]
and
\[
J_i^R \rightharpoonup J_i \textnormal{ weakly in } L^2(\textnormal{supp}\phi;\mathbb{R}^3)
\]
we may pass to the limit to obtain \eqref{weakmass}, or 
\[
\del_t c_i + \diver J_i = 0
\]
in distributions.

\subsection*{Passage to the limit in the linear system} 

Let $\psi\in C_c^\infty(\mathbb{R}^3\times[0,T);\mathbb{R}^3)$ be a test function, and choose $K>0$ such that $\textnormal{supp}\psi \subset B_K\times[0,T)$. For $R>K$, the weak formulation of the linear system reads
\[
\int_0^T\int_{\mathbb{R}^3} c_i^R \diver \psi \,dxdt = \int_0^T\int_{\mathbb{R}^3} \sum_{j\not=i} \frac{c_j^RJ_i^R-c_i^RJ_j^R}{D_{ij}} \cdot \psi \,dxdt.
\]
The left-hand side passes to the limit because of the strong convergence of $c_i^R$. For the right-hand side, the strong convergence of $c_i^R$ and the weak convergence of $J_i^R$, imply that
\[
c_i^RJ_j^R \rightharpoonup c_iJ_j \textnormal{ weakly in } L^2(\textnormal{supp}\psi;\mathbb{R}^3)
\]
and passing to the limit, we obtain
\[
\int_0^T\int_{\mathbb{R}^3} c_i \diver \psi \,dxdt = \int_0^T\int_{\mathbb{R}^3} \sum_{j\not=i} \frac{c_jJ_i-c_iJ_j}{D_{ij}} \cdot \psi \,dxdt,
\]
or \eqref{linearsystem} in distributions.

\subsection*{Passage to the limit in the entropy inequality}

It remains to pass to the limit in the entropy inequality. Since
\[
c^R \to c \textnormal{ a.e. locally}
\]
and $\tilde c^R = \bar c$ outside $B_R$, we have that
\[
\tilde c^R \to c \textnormal{ a.e. in } \mathbb{R}^3\times(0,T).
\]
By Fatou's Lemma
\[
H[c(t)|\bar c] \leq \liminf_{R\to\infty} H_{B_R}[c^R(t)|\bar c].
\]

For the gradient term, we have that 
\[
\sqrt{c_i^R} \to \sqrt{c_i} \textnormal{ strongly in } L^2_{\rm loc}(\mathbb{R}^3\times(0,T))
\]
due to the inequality
\[
|\sqrt{u}-\sqrt{v}|^2 \leq |u-v| \textnormal{ for all } u,v\in[0,1]
\]
and the fact that $c_i^R\to c_i$ strongly. Moreover
\[
\nabla\sqrt{c_i^R} \rightharpoonup \nabla\sqrt{c_i} \textnormal{ weakly in } L^2_{\rm loc}
\]
thus for every $K>0$
\[
\int_0^t\int_{B_K} |\nabla\sqrt{c_i}|^2 \,dxds \leq \liminf_{R\to\infty} \int_0^t\int_{B_K} |\nabla\sqrt{c_i^R}|^2 \,dxds.
\]
Letting $K\to\infty$, by the Monotone Convergence Theorem we find
\[
\int_0^t\int_{\mathbb{R}^3} |\nabla\sqrt{c_i}|^2 \,dxds \leq \liminf_{R\to\infty} \int_0^t\int_{B_R} |\nabla\sqrt{c_i^R}|^2 \,dxds.
\]
Combining the uniform entropy inequality \eqref{entrineq}, with the control of the gradient by the dissipation \eqref{grad}, we obtain
\[
H_{B_R}[c^R(t)|\bar c] + C\int_0^t\int_{B_R} \sum_i|\nabla\sqrt{c_i^R}|^2 \,dxds \leq H[c^0|\bar c]
\]
and passing to the limit as shown above, we arrive at
\[
H[c(t)|\bar c] + C\int_0^t\int_{\mathbb{R}^3} \sum_i|\nabla\sqrt{c_i}|^2 \,dxds \leq H[c^0|\bar c]
\]
for almost every time $t\in(0,T)$.

Finally, by weak lower semicontinuity
\[
\int_0^T\int_{\mathbb{R}^3} |J_i|^2 \,dxdt \leq \liminf_{R\to\infty} \int_0^T\int_{B_R} |J_i^R|^2 \,dxdt \leq C H[c^0|\bar c]
\]
and thus
\[
J_i \in L^2(\mathbb{R}^3\times(0,T);\mathbb{R}^3).
\]


\end{document}